\documentclass[papera4,12pt]{article}
\usepackage[english]{babel}
\usepackage[letterpaper,top=2cm,bottom=2cm,left=3cm,right=3cm,marginparwidth=1.75cm]{geometry}
\usepackage{amsmath}
\usepackage{graphicx}
\usepackage[english]{babel}
\usepackage[colorlinks=true, allcolors=red]{hyperref }
\usepackage{amsthm}
\usepackage{amsfonts}
\usepackage{hyperref}
\newtheorem{theorem}{Theorem}[section]
\newtheorem{lemma}{Lemma}[section]
\numberwithin{equation}{section}
\begin{document}
\begin{center}
{\bf {\Large ON THE DIOPHANTAINE EQUATIONS

$J_n +J_m =F_a$ 

$F_n +F_m =J_a$ 

}}
\vspace{8mm}

{\Large \bf Seif Tarek$^1$, Ahmed Gaber$^2$ and  M. Anwar$^3$}
\vspace{3mm}
$^1$$^2$$^3$ 

Department of Mathematics, University of Ain Shams \\ 
Faculty of science, Egypt \\
e-mails:\href{mailto:starek98@sci.asu.edu.eg} {\url{starek98@sci.asu.edu.eg}}$^1$,\\ \href{mailto:a.gaber@sci.asu.edu.eg}{\url{a.gaber@sci.asu.edu.eg}}$^2$,\\  
\href{mailto:mohmmed.anwar@hotmail.com}{\url{mohmmed.anwar@hotmail.com}},\hspace{0.5cm}\href{mailto:mohamedanwar@sci.asu.edu.eg}{\url{mohamedanwar@sci.asu.edu.eg}}$^3$  
\end{center}
\vspace{10mm}

\begin{abstract}
Let $\left\lbrace F_{k}\right\rbrace_{k\geq0}$ be the  Fibonacci sequence defined by  $F_{k}=F_{F-1}+F_{k-2}$ for all $ n\geq2$ with initials  $F_{0}=0\; F_{1}=1$. Let $\left\lbrace J_{n}\right\rbrace_{n\geq0}$ be the Jacobsthal sequence defined by $J_n=2J_{n-1}+J_{n-2}$ for all $ n\geq2$ with initials $J_0=0$, $J_1=1$. In this paper we find all the solutions of the two Diophantine equations $J_n +J_m =F_a$ ,$F_n +F_m =J_a$ in the non-negative integer variables (n,m,a),i.e we determine all Fibonacci numbers which are sum of two Jacobsthal numbers, and also determine all Jacobsthal numbers which are sum of two Fibonacci numbers.
\end{abstract}
\begin{flushleft}
\textbf{2020 Mathematics Subject Classification}: 11B83.\\
\end{flushleft}
\begin{flushleft}
\textbf{Keywords}: 

\end{flushleft}
\section{Introduction}

the Fibonacci numbers, commonly denoted by $F_n$, form sequence, the Fibonacci sequence, where every number in it is the sum of the two preceding ones , it starts with $F_0 =0$ , $F_1 =1$ , so we can define the Fibonacci sequence as follows,\\

\textbf{Definition}:the Fibonacci sequence $(F_n)_{n\geq 0}$ given by $F_0=0$,$F_1=1$ and 

\begin{equation*}
F_{n+2}=F_{n+1}+F_{n} \;\; for \; all\; n\geq0 . 
\end{equation*}

The first 15 Fibonacci numbers are  0,1,1,2,3,5,8,13,21,34,55,89,144,233,377,610,
the Fibonacci numbers were first mentioned in indian mathematics ,as early as 200 BC by pingala , they  are named after the Italian mathematician leonardo of pisa , later they known as Fibonacci, who introduced the sequence  to western European Mathematics .
the Fibonacci numbers appear unexpectedly often in Mathematics, the applications of Fibonacci numbers include computer algorithms,  biological setting,and many other fields.
also Fibonacci numbers are strongly related to the golden ratio , as the ratio of two consecutive Fibonacci numbers tends to the golden ratio as n increases.\\

Jacobsthal polynomials $(J_n(x))_{n\geq0}$ were first studied by E.E.Jacobsthal around 1919 , this sequence can be defined by the recurrence $J_0(x)=0$, $J_1(x)=1$ and $J_n(x)=J_{n-1}(x)+J_{n-2}(x)$ for all $n\geq 0$, we can see that, $J_n(1)=F_n$, the Fibonacci sequence.the Jacobsthal sequence $(J_n)_{n\geq0}$ can be defined as $J_n=J_n(2)$,also we can defined the Jacobsthal sequence by the recurrence   
$J_0=0$, $J_1=1$ and $J_{n+2}(x)=J_{n+1}(x)+2 J_{n}(x)$ for all $n\geq 0$.the first 10 Jacobsthal numbers are 0,1,1,3,5,11,21,43,85,171,341,... . the first apperance of the Jacobsthal sequence, Horadam was the first who considered such a sequence in detail in his seminal paper, Horadam's work motivates a lot of researchs.
%for some recent works.

The aim of this work is to determine all Fibonacci numbers which are sum of two Jacobsthal numbers, and also determine all Jacobsthal numbers which are sum of two Fibonacci numbers,i.e we determine all non negative integer solutions (n,m,a) of the following two Diophantine equations 

\begin{equation}\label{eqn1}
J_n +J_m =F_a   
\end{equation}

\begin{equation}
F_n +F_m =J_a
\end{equation}

and our results are the following.  

\begin{theorem}
  The only solutions of the Diophantaine equation $J_n +J_m =F_a$ in the nonnegative integers n,m and a with $n \geq m$ are given by \\

\begin{equation*}
(n,m,a)\in   \left\{ \begin{array}{cc}
        & (0,0,0),\hspace{0.5cm}(1,0,1),\hspace{0.5cm}(1,0,2),\hspace{0.5cm}(2,0,2),\hspace{0.5cm}(1,1,3), \\
        & (2,1,3),\hspace{0.5cm}(2,2,3),\hspace{0.5cm}(3,0,4),\hspace{0.5cm}(4,0,5),\hspace{0.5cm}(4,3,6),\\
        & (6,0,8)
   \end{array}
   \right\}.
\end{equation*}
\end{theorem}

\begin{theorem}
The only solutions of the Diophantine equation $F_n +F_n =J_a$ in the non negative integers n,m and a with $n \geq m$ are given by

\begin{equation*}
  (n,m,a)\in \left\{ \begin{array}{cc}
        & (0,0,0),\hspace{0.5cm}(1,0,1),\hspace{0.5cm}(1,0,2),\hspace{0.5cm}(2,0,1),\hspace{0.5cm}(2,0,2), \\
        & (4,0,3),\hspace{0.5cm}(5,0,4),\hspace{0.5cm}(8,0,6),\hspace{0.5cm}(3,1,3),\hspace{0.5cm}(3,2,3),\\
        & (4,3,4),\hspace{0.5cm}(6,4,5),\hspace{0.5cm}(7,6,6),\hspace{0.5cm}(8,0,6)
   \end{array}
   \right\}.
\end{equation*}

\end{theorem}
The  strategy for providing Theorem 1.1 is as follows :First,we find upper bound for the nonnegative integer a as a function of n, second, we rearrange Equation (1.1) in suitable ways in order to get two different linear forms in logarithms of algebraic numbers which are both nonzero and small.Next we use twice a lower bound on such nonzero linear forms in logarithms due to Matveev to find an absolute upper bound for n,hence,an absolute upper bound for a and m,then we reduce this upper bound using standard facts about continued fractions.The strategy for providing Theorem (1.2) is quite similar.

\section{preliminary Results}

\subsection{The Binet formula for the Fibonacci numbers}
Like every sequence defined by a linear recurrence with constant coefficients the Fibonacci numbers have a closed form expression known as Binet's formula  see (\cite{1})

\begin{equation*}
F_n=\frac{\alpha^n-\beta^n}{\sqrt{5}} \hspace{0.25cm} \text{holds for all} \; n\geq 0
\end{equation*}

where 

\begin{equation*}
\alpha=\frac{1+\sqrt{5}}{2}\; \hspace{0.25cm} \;\beta= \frac{1-\sqrt{5}}{2} .
\end{equation*}

also by induction  we can prove easily that the inequality 
\begin{equation}
    \alpha^{n-2} \leq F_n \leq \alpha^{n-1} 
\end{equation}
holds for all $n\geq1$.

we will need this relations in our main theorems .

\subsection{Binet formula for the Jacobsthal numbers }
Like the Fibonacci numbers,Jacobsthal numbers have a closed form (see\cite{1}). Since
\begin{equation*}
J_n(x)=\frac{u^n-v^n}{u-v} \hspace{0.25cm}\text{holds for all}\; n\geq 0 []
\end{equation*}

where 
\begin{equation*}
u=u(x)=\frac{1+\sqrt{4x+1}}{2}\; ,\; v=v(x)=\frac{1-\sqrt{4x+1}}{2} \;.\\
\end{equation*}

therefore 

\begin{equation*}
J_n=\frac{2^n-(-1)^n}{3} \hspace{0.25cm} \text{holds for all} \hspace{0.25cm} n\geq 0 .
\end{equation*}
also we can prove the next inequality by induction   
\begin{equation}
2^{n-2}<J_n<2^{n-1} \hspace{0.25cm} \text{for all }\hspace{0.25cm} n\geq 3
\end{equation}

now we recall some basic notions from algebraic number theory .

\subsection{Linear forms in logarithms}

let $\gamma$ be an algebraic number of degree d and its minimal polynomial is 
\begin{equation*}
a_0 x^d +a_1 x^{d-1} + \cdots + a_d = a_0 \prod_{i=1}^{d}(x-\gamma^{(i)})
r\end{equation*}

 where ${a_i}$'s are relativity prime integers and $a_0>0$ and $\gamma^{(i)}$'s are the conjugates of  $\gamma$, then the logarithmic height of $\gamma$ is defined as

\begin{equation*}
h(\gamma)= \frac{1}{d} \big(\log {a_0} + \sum_{i=1}^{d}\log(max\{|\gamma^{(i)}|,1\})\big)
\end{equation*}
we can see that if $\gamma = p/q$ is rational number with $gcd(p,q)=1$,then $h(\gamma)=\log max\{|p|,q\}$.

the  following  properties of the logarithmic height will be used in the main theorems without special reference \\

(i) $h(z+s)\leq h(z)+h(s)+\log 2$\\

(ii) $h(z\: s^{\pm 1})\leq h(z)+h(s)$\\

(iii) $h(z^s)\leq |s|\: h(z)$

with the previous notions Matveev (see\cite{3}) proved the next theorem .

\begin{theorem}
    
(Matveev’s theorem) Assume that $\gamma_1,\cdots , \gamma_t$ are positive real alge-
braic numbers in a real algebraic number field K of degree D, $b_1,\cdots , b_t$ are rational
integers, and

\begin{equation*}
\Gamma =\gamma^{b_1} \cdots \gamma^{b_t} -1 \neq 0
\end{equation*}
 then  
 \begin{equation}
 \log |\Gamma| > -1.4 \times 30^{t+3} \times t^{4.5} \times D^2 \: (1+\log D) (1+\log B)\: A_1\:A_2\:A_3
 \end{equation}
where $B\geq max\{|b_1|\cdots|b_t|\}$ , and
$A_i \geq max\{D \: h(\gamma_i) , |\log\gamma_i| , 0.16\}$\\ \hspace{0.5cm} for all $i=1,\cdots , t$.
\end{theorem}

Dujella and peth{\"o} in \cite{2} had  a version of the reduction method based on the Baker-Davenport lemma in . the next lemma is a key tool for reducing the upper bounds on the variables of equations (0,1)and(0.2).   

\begin{lemma}
Suppose that M is a positive integer, and p/q is a convergent of the continued fraction of the irrational number $\gamma$ ,such that $q>6\: M$ , and A,B,$\mu$ are some real numbers with $A>0$ and $B>1$ , let $\epsilon =||\mu \:q|| -M\:||\gamma \: q|| $ ,where $||.||$ denotes the distance from the nearest integer. if $\epsilon >0$ then there exist no solution to the inequality 

\begin{equation*}
0<|u \: \gamma - v + \mu|<A\:B^{-w}    
\end{equation*}

in the positive integers u , v  and w with  
\begin{equation*}
u\geq M \; and \; w\geq \frac{\log(Aq / \epsilon)}{\log B}
\end{equation*}

\end{lemma}

now we prove our main theorems .

\section{Proof of Theorem 1.1}
\begin{proof}

Because of the symmetry of equation (1.1) we can assume that $n \geq m$, by using Mathematica we found that the only solutions of equation (1.1) in the range $m\leq n \leq200$, are the mentioned solutions in Theorem 1.1. We now assume that $n>200$, and we divide the proof into three steps.   

\subsection{Finding relations between a and n.} 

we begin by proving that    
\begin{equation*}
F_n < J_n \hspace{0.5cm} \text{holds for all} \hspace{0.5cm} n > 200.
\end{equation*}
By induction we can prove the following inequality 
\begin{equation*}
\alpha^{n-1} \leq 2^{n-2}\; \text{holds for all}\; n > 200 .
\end{equation*}
combining the last relation with inequality (2.1) and inequality (2.2) we get that 
\begin{equation*}
F_n \leq \alpha^{n-1} \leq 2^{n-m} < J_n  .
\end{equation*}

holds for all $n > 200$,if $n\geq a$ then $J_n > F_a$,which is false since $J_n +J_m =F_a$.therefore 

\begin{equation}
a > n 
\end{equation}
combining (1.1) with the lift hand side of (2.1) and the right hand side of (2.2) we get that 
\begin{equation*}
\alpha^{a-2} \leq F_a = J_n + J_m \leq 2J_n < 2 \times 2^{n-1} = 2^n 
\end{equation*}

Taking logarithms in the last inequality and rearranging we get 

\begin{equation*}
a < \frac{\log2}{\log\alpha}  n + 2  
\end{equation*}

since\;$n\geq 200$ ,it follows that  
\begin{equation}
n < a < 1.6n
\end{equation}

\subsection{Finding upper bound on both n and m. }
now we rewrite equation (1.1) as fallows
\begin{equation*}
\frac{2^n}{3} - \frac{\alpha^a}{\sqrt{5}} = -J_m - \frac{\beta^a}{\sqrt{5}} + \frac{(-1)^n}{3}
\end{equation*}
taking the absolute values to the above relation we get that
\begin{equation*}
|\frac{2^n}{3} - \frac{\alpha^a}{\sqrt{5}}| = |-J_m - \frac{\beta^a}{\sqrt{5}} + \frac{(-1)^n}{3}| \leq J_m+ \frac{|\beta|^a}{\sqrt{5}} + \frac{1}{3}
\end{equation*}
Since $J_m \leq 2^m$ , for all $m\geq 0$ , and $\frac{|\beta|^a}{\sqrt{5}} < \frac{1}{2}$ , for all $a > 200 $ , it follows that 
\begin{equation*}
|\frac{2^n}{3} - \frac{\alpha^a}{\sqrt{5}}| < 2^m + \frac{5}{6} 
\end{equation*}

Dividing both sides of the above expression by $2^n/3$  we get 
\begin{equation}
|1-2^{-n} \; \alpha^a \; \frac{3}{\sqrt{5}}| < \frac{5}{2^{n-m}}
\end{equation}
In order to apply matveev's theorem we first prove that the left hand side of (3.3) is not zero,If this were zero,we would then get that $ 9 \alpha^{2a}$ is an integer which is impossible . We take The parameter $t=3$ ,and $\gamma_1=2$, $\gamma_2=\alpha$ and $\gamma_3=3/\sqrt{5}$,We also take $b_1=-n$,$b_2=a$ and $b_3=1$. since max$\{|-n|,|a|,|1|\}=a$,then\; $B=a$ ,and the smallest field contains $\gamma_1 ,   \gamma_2 and \; \gamma_3$ is $\mathbb{Q}$$(\sqrt{5})$, so\;$D=2$. $h(\gamma_1)=\log 2$ therefore we can choose $A1=1.4>D h(\gamma_1)$ , $h(\gamma_2)=(\log \alpha)/2$ , therefore we can take $A2=0.5>D h(\gamma_2)$ , $h(\gamma_3)=h(\frac{3}{\sqrt{5}})=\frac{1}{2}(\log5 + 2\log \frac{3}{\sqrt{5}}) =\log 3$ , therefore we can take $A3=2.2>D h(\gamma_3)$. Applying Matveev's theorem to (3.3), we get 
\begin{equation*}
    \log|1-2^{-n}\; \alpha^a\; \frac{3}{\sqrt{5}}| > -C(1+\log a) \times1.4 \times 0.5 \times 2.2.
\end{equation*}
where $C= 1.4 \times 30^6 \times 3^{4.5} \times 2^2\times (1+log2) < 9.7\times 10^{11}$.
Using the fact that $2\log a > 1+\log a$ for every $a\geq 3$,
 And taking logarithms in inequality (1.5) and comparing the resulting inequality with the above inequality ,we get
\begin{equation*}
\log5 -(n-m)\log2 > -2.99 \times 10^{12} \times \log a 
\end{equation*}
The last inequality implies that 
\begin{equation}
(n-m) \log2 < 3 \times 10^{12}\times \log a. 
\end{equation}
Now we consider another linear form of (1.1). To this end, We rewrite equation(1.1) as fallows \begin{equation*}
\frac{2^n}{3}+\frac{2^m}{3}-\frac{\alpha^a}{\sqrt{5}}=\frac{-\beta^a}{\sqrt{5}}+\frac{(-1)^n}{3}+\frac{(-1)^m}{3}.
\end{equation*}
taking the absolute values of the above inequality,we get 

\begin{equation*}
    |\frac{2^n}{3}+\frac{2^m}{3}-\frac{\alpha^a}{\sqrt{5}}|=|\frac{-\beta^a}{\sqrt{5}}+\frac{(-1)^n}{3}+\frac{(-1)^m}{3}| < \frac{1}{2} +\frac{1}{3} + \frac{1}{3} =\frac{7}{6}
\end{equation*}
we used that 
$\frac{|\beta|^a}{\sqrt{5}}<\frac{1}{2}$ for every $a>200$,rearranging the above inequality we get

\begin{equation*}
|\frac{2^n}{3} (1+2^{m-n}) - \frac{\alpha^a}{\sqrt{5}}| < \frac{7}{6}
\end{equation*}
Dividing by the first term of the left hand side of the above inequality,we get
\begin{equation}
|1-\alpha^a\; 2^{-n}\; \frac{3}{\sqrt{5}} (1+2^{m-n})^{-1}| < \frac{7}{2}\times \frac{1}{2^n}  
\end{equation}
Now in order to apply matveev's theorem,we prove that the left hand side of (1.7) is not equals to zero , if it equals zero then we would get that 
\begin{equation*}
(2^n + 2^m) \sqrt{5} = 3\alpha^a
\end{equation*}
Conjugating both sides of the above inequality,we get 
\begin{equation*}
-(2^n + 2^m) \sqrt{5} = 3\beta^a 
\end{equation*}
which implies that  $|\alpha^a| = |\beta^a|$,therefore we get that
\begin{equation*}
\alpha^a = |\alpha^a| = |\beta^a| < |\beta|^a < 1
\end{equation*}
 which is impossible for any positive integer a. Applying matveev's theorem to (3.5),by Taking The parameter $t=3$ , $\gamma_1=2, \gamma_2=\alpha , \gamma_3= \frac{3}{\sqrt{5}} (1+2^{m-n})^{-1} $, we also take $b_1 = -n , b_2 = a , b_3 = 1 $.
then we can take  $B=a$ ,and the smallest field contains $\gamma_1 ,   \gamma_2\;and\; \gamma_3$ \;is $\mathbb{Q}$$(\sqrt{5})$\; so\; $D=2$. since $h(\gamma_1)=\log 2$, we can choose $A_1=1.4>Dh(\gamma_1)$ , also since $h(\gamma_2)=\frac{\log\alpha}{2}$ , we can take $A_2=0.5>Dh(\gamma_2)$ , now we estimate $ h (\gamma_3)$.

\begin{eqnarray*}
 h (\gamma_3)=h(\frac{3}{\sqrt{5}} (1+2^{m-n})^{-1})&\leq &h(\frac{3}{\sqrt{5}}) + h(2^{m-n}) + \log2 \\ &=&\log3 + |m-n|\: h(2)+\log2 
\\&=&\log6 + (n-m) \log2
\end{eqnarray*} 

Next, notice that 
\begin{equation*}
|\log\gamma3| = |\frac{3}{\sqrt{5}} (1+2^{m-n})^{-1}| < |\log  \frac{3}{\sqrt{5}}| + |\log(1+2^{m-n})^{-1}|<\log \frac{3}{\sqrt{5}} +\log2< 1
\end{equation*}
hence, we can take
\begin{equation*}
A_3=4 + 2(n-m) \log2 >max\{2h(\gamma_3) , |\log\gamma_3| , 0.16\}
\end{equation*}
Applying Matveev's theorem to (3.5), we get  
\begin{equation*}
\log |1-\alpha^a\; 2^{-n}\; \frac{3}{\sqrt{5}} (1+2^{m-n})^{-1}| > -C (1+log a)\times1.4 \times 0.5 \times (4 + 2(n-m) \log2) 
\end{equation*}
where $C=1.4 \times 30^6 \times 3^{4.5} \times 2^2\times (1+log2)<9.7\times10^{11}$. using the fact that\; $2 \log a > 1+\log a$ \;for every\; $a\geq 3$,we get 
\begin{equation}
\log |1-\alpha^a\; 2^{-n}\; \frac{3}{\sqrt{5}} (1+2^{m-n})^{-1}| > -1.36 \times 10^{12} \times (\log a)\times (4 + 2(n-m) \log2)
\end{equation}
taking logarithms in inequality (3.5) and comparing it with(3.6), we get 

\begin{equation}
\log\frac{7}{2} - a\log2 > -1.36 \times 10^{12}\times (\log a) \times(4 + 2(n-m) \log2)
\end{equation}
which implies that 
\begin{equation}
a\log2 <1.37 \times 10^{12}\times (\log a) \times(4 + 2(n-m) \log2)
\end{equation}
substituting  by (3.4) in (3.8), we get 
\begin{equation}
a \log2 < 1.37 \times10^{12}\times (\log a) (4 +6 \times 10^{12} \log a ) 
\end{equation}

solving inequality (1.13) by  Mathematica, we get 
\begin{equation}
a < 10^{29} 
\end{equation}
And since $n<a$, we get
 \begin{equation}
 n < 10^{29}
\end{equation}
\subsection{Reducing the bound on n.}
Now after finding an upper bound on n, we reduce it to a size that can be easily handled, we will use Lemma 1 several times to achieve that. First , let   

\begin{equation}
Z = a \log\alpha - n \log2   + \log\frac{3}{\sqrt{5}}
\end{equation}
Note that (3.3) can be rewritten as 
\begin{equation}
|e^z - 1| < \frac{5}{2^{n-m}}
\end{equation}
Since we proved that the left hand side of (3.3) is not zero,we get that $Z \neq 0$.   
If $Z > 0$ then
\begin{equation}
0 < Z < e^z - 1 =|e^z - 1| < \frac{5}{2^{n-m}}
\end{equation}
If $Z < 0$, we suppose that $n-m \geq 20$, then we get that 
\begin{equation}
e^z - 1< \frac{1}{2}
\end{equation}
In the above inequality we used the inequality $\frac{5}{2^{n-m}} < \frac{1}{2}$, where $n-m \geq20$. Inequality (3.15) implies that $e^{-z} < \frac{1}{2}$, or $e^{|z|} < 2$, then we get  
\begin{equation}
0 < |Z| < e^{|z|} - 1 = e^{|z|} (1-e^z) = e^{|z|} |e^z -1| < \frac{10}{\alpha^{n-m}}
\end{equation}
 in both cases $(z > 0 , z < 0)$, (3.16) is true , replacing Z in the above inequality by its formula (3.12), and dividing by $log2$, we get 
\begin{equation} 
0 < |a\frac{log \alpha}{log2} - n + \frac{log\frac{3}{\sqrt{5}}}{log2}| < \frac{15}{2^{n-m}}
\end{equation}
We put 
\begin{equation*}
\gamma = \frac{\log \alpha}{\log 2},\hspace{1cm} \mu = \frac{\log\frac{3}{\sqrt{5}}}{\log2},\hspace{1cm} A = 15,\hspace{1cm} B = 2,\hspace{0.5cm}\text{and}\hspace{0.5cm} W = n-m.
\end{equation*}
Clearly $\gamma$ is irrational , and from (3.10) we can take $M = 10^{29}$, as an upper bound on a. For $q_{(69)} = 20721505928824926197089563175427> 6M $, Where q is a denominator of a convergent of the continued fraction of $\gamma$, Mathematica found that,  $0.333233182722303< \epsilon=||\mu \:q|| -M\:||\gamma \: q||<0.333233182722304$ , applying Lemma 2.1 to inequality (3.17), we get 
\begin{equation}
n-m<\frac{log(Aq/\epsilon)}{log B} < 109.53
\end{equation}
Which implies that 
\begin{equation}
n-m \in[1,109]    
\end{equation}

Inserting this upper bound for  $n-m $ into (3.8), we get 
\begin{equation}
a < 1.2 \times10^{16} 
\end{equation}
and so 
\begin{equation}
n < 1.2 \times 10^{16} .
\end{equation}
Now we work on (3.5) in order to find a better upper bound on n.

taking 
\begin{equation}
s = a \log\alpha - n \log2 + \log\frac{3}{\sqrt{5}} (1+2^{m-n})^{-1} 
\end{equation}
Inequality (1.9) implies that 
\begin{equation}
|1-e^s|<\frac{7}{2}\times \frac{1}{2^n}
\end{equation}
Note that $s \neq 0$, since$|1-e^s|\neq 0$.    

If $s > 0$, then 
\begin{equation}
0 < s < e^s - 1 =|e^s - 1| < \frac{7}{2} \times \frac{1}{2^n}
\end{equation}
And if $s < 0$, we get  
\begin{equation}
(1-e^s)=|e^s-1|  < \frac{7}{2} \times\frac{1}{2^n}< \frac{1}{2} \;\: \; \text{for all}\; n > 200.
\end{equation}
The last inequality implies that $e^{-s} < 2$, or \;$e^{|s|} < 2$\; which implies that 
\begin{equation}
0 < |s| < e^{|s|} - 1 = e^{|s|} (1-e^s) = e^{|s|} |e^s -1| < \frac{7}{2^n}
\end{equation}
In both cases $(s > 0 , s < 0)$ inequality (3.26) is true. Replacing the value of s in the above inequality by its formula (3.22) and dividing by $\log2$, we get
\begin{equation}
0 < |a\frac{\log \alpha}{\log2} - n + \frac{\log{\frac{3}{\sqrt{5}} (1+2^{m-n})^{-1}}}{\log2}| < \frac{11}{2^n}
\end{equation}
We put 
\begin{equation*}
\gamma = \frac{\log \alpha}{\log2},\hspace{1cm} \mu_{(n-m)} =\frac{\log{\frac{3}{\sqrt{5}} (1+2^{-(n-m)})^{-1}}}{\log2},\hspace{1cm} A = 11,\hspace{1cm}B = 2,\hspace{0.5cm}\text{and}\hspace{0.5cm} W = n.
\end{equation*}
Using inequality (3.20) we can Take  $M=1.2 \times 10^{16}$ as an upper bound on a, taking
 $q_{(82)}=1234165504911193651820557190855668171489 > 6M $, with the help of Mathematica,we find that for all choices $n-m \in[1,109]$ the smallest epsilon is $0.00663531736488705<\epsilon=||\mu_{(66)} \:q_{(82)}|| -M\:||\gamma \: q_{(82)}|| < 0.00663531736488707$, which belongs to $n-m=66$, so applying lemma 2.1 to inequality (3.27) ,we get
\begin{equation}
n<\frac{log(Aq/\epsilon)}{log B} < 140,56
\end{equation}
In the last inequality we used the smallest epsilon because it produces the largest upper bound on n,(3.28) gives a contradiction to our assumption that $n \leq 200$, which shows that there are no solutions to (1.1) when $n > 200$. \newpage
\end{proof}

\section{Proof of Theorem 1.2}

\begin{proof}
First of all we can assume that $n\geq m$, because of the symmetry of equation (1.2). By a search in Mathmatica we find that the only solutions of equation (1.2) in the range $m\leq n \leq 200$ are the mentioned solutions in Theorem 1.2. So we assume now that $n>200$ and try to find if there are other solutions. As in the proof of Theorem 1.1 we divide the proof into three steps.    

\subsection{Finding relations between a and n.}
First, Theorem 1.1 tells Us that the only solutions of the equation $F_n = J_a$ are $(F_0=J_0,F_1=J_1,F_2=J_1,F_2=J_2,F_4=J_3,F_5=J_4,F_8=J_6)$, which implies that for $n>200$ the nonnegative integer m must be greater than zero, i.e $m\geq 1$. now we begin by proving that 
\begin{equation}
   2F_n<J_n \hspace{1cm} \text{for all}\hspace{0.5cm} n>200
\end{equation}
by induction it is easy to prove that 
\begin{equation*}
2\alpha^{n-1}<2^{n-2} \hspace{1cm}\text{for all}\hspace{0.5cm} n>200.
\end{equation*}
Combining the last relation with (2.1) and (2.2), we get 
\begin{equation*}
2F_n< 2\alpha^{n-1} <2^{n-2}<J_n.
\end{equation*}
Inequality (4.1) implies that 
\begin{equation}
    n>a.
\end{equation}  
to see this let us assume that $n\geq a$ then $F_n+F_m\leq 2F_n < J_n\leq J_a $, which contradicts equation (1.2). 

\subsection{Finding upper bound on n}
now we rewrite equation (1.2) as 

\begin{equation*} 
\frac{\alpha^{n}}{\sqrt{5}}-\frac{2^a}{3}=-F_m-\frac{(-1)^a}{3}+\frac{\beta^n}{\sqrt{5}}
\end{equation*}
Take the absolute values in the above relation,  we get
\begin{equation*}
|\frac{\alpha^{n}}{\sqrt{5}}-\frac{2^a}{3}|=|-F_m-\frac{(-1)^a}{3}+\frac{\beta^n}{\sqrt{5}}|\leq F_m+\frac{|\beta|^n}{\sqrt{5}}+\frac{1}{3}<\alpha^m+\frac{5}{6}
\end{equation*}
In the last step we used the fact that $\frac{|\beta|^n}{\sqrt{5}}<\frac{1}{2}$ for every $n>200$, and $F_m<\alpha^m$, dividing the above relation by $\frac{\alpha^{n}}{\sqrt{5}}$,we get
\begin{equation}
|1-2^a\; \alpha^{-n}\; \frac{\sqrt{5}}{3}|<\frac{4}{\alpha^{n-m}}
 \end{equation}
In order to use Matveev's Theorem we prove that the left hand side of (2.3) is not equal zero, 
 If it equals zero, then $\sqrt{5}\: 2^a=3\alpha^a$, which implies that $9\alpha^{2n}$ is an integer which is impossible. now we apply matveev's Theorem, taking the parameter $t=3$ , $\gamma_1=2, \gamma_2=\alpha , \gamma_3=\frac{\sqrt{5}}{3}$, we also take, $b_1=a,b_2=-n,b_3=1$. Since $B=max\{|b_1|,|b_2|,|b_3|\}$, $B=n$ and the smallest field contains $\gamma_1 ,   \gamma_2 , \gamma_3$ , is $\mathbb{Q}$$(\sqrt{5})$\; so\; $D=2$. Since $h(\gamma_1)=\log 2$,  we can choose $A1=1.4>Dh(\gamma_1)$, also since  $h(\gamma_2)=\frac{\log\alpha}{2}$, we can take $A2=0.5>Dh(\gamma_2)$, similarly $h(\gamma_3)=h(\frac{\sqrt{5}}{3})=\log 3$, so we can take $A3=2.2>Dh(\gamma_3)$, applying Matveev's theorem to inequality (4.3), we get 
\begin{equation*}
\log|1-2^a\; \alpha^{-n}\; \frac{\sqrt{5}}{3}|>-1.4\times30^6 \times 3^{4.5}\times 2^2 (1+\log2) \times (1+\log n) \times1.4 \times0.5 \times 2.2
\end{equation*}
Using the fact that $2\log n > 1+\log n$ for every $n\geq 3$,
and taking logarithms in inequality (4.3) and comparing the resulting inequality with the above inequality, we get 
\begin{equation*}
\log4 - \log\alpha^{n-m} > -2.987 \times10^{12} \times \log n .
\end{equation*}
the last inequality implies that 
\begin{equation}
(n-m )\log\alpha < 3 \times10^{12} \log n
\end{equation}

Now we back again to (1.2) and try to find a second linear form , so we rewrite equation (1.2) as follows
 \begin{equation*}
\frac{\alpha^n}{\sqrt{5}}+\frac{\alpha^m}{\sqrt{5}}-\frac{2^a}{3}=\frac{\beta^n}{5}+\frac{\beta^m}{5}-\frac{(-1)^a}{3}
.\end{equation*} 
Taking the absolute values in the above relation, we get
 \begin{equation*}
|\frac{\alpha^n}{\sqrt{5}} (1+\alpha^{m-n}) - \frac{2^a}{3}| \leq \frac{|\beta|^{n}+|\beta|^{m}}{\sqrt{5}}+\frac{1}{3}<\frac{2}{3}
\end{equation*}

We used in the above relation the fact that $\frac{|\beta|^{n}+|\beta|^{m}}{\sqrt{5}}<\frac{1}{3} ,\forall n\geq 5 , m\geq 1$. Dividing the above inequality by the first term in its left hand side, we get 
\begin{equation}
|1-2^a \;\alpha^{-n}\; \frac{\sqrt{5}}{3} (1+\alpha^{m-n})^{-1}| < \frac{2}{\alpha^n}
\end{equation}

Now we prove that the left hand side of (4.5) is not equals zero. If it equals zero, then $3\alpha^n (1+\alpha^{m-n}) = \sqrt{5} 2^a$ , which implies that 

\begin{equation*}
3(\alpha^n + \alpha^m) = \sqrt{5} 2^a
\end{equation*}

conjugating both sides, we get  

\begin{equation*}
3(\beta^n + \beta^m) = -\sqrt{5} 2^a
\end{equation*}
Combining the last two relations, we get 
\begin{equation*}
\alpha^n < \alpha^n + \alpha^m = |\alpha^n + \alpha^m| = |\beta^n + \beta^m | < |\beta|^n + |\beta|^m < 1
\end{equation*}

which is impossible for any positive integer n, taking $\gamma_1 = 2 , \gamma_2 = \alpha , \gamma_3= \frac{\sqrt{5}}{3} (1+\alpha^{m-n})^{-1} , b_1 = a , b_2 = -n , b_3 = 1$ and so  we take the parameter $t=3$, and Clearly the smallest field contains $(\gamma_1 , \gamma_2 , \gamma_3)$\;is $\mathbb{Q}$$(\sqrt{5})$\; so\; $D=2$, also $B=max\{|b_1|,|b_2|,|b_3|\}=n$. As before $A1=1.4$, $A2=0.5$, now we estimate $h(\gamma3)$.
\begin{eqnarray*}
h(\frac{\sqrt{5}}{3} (1+\alpha^{m-n})^{-1})\leq h(\frac{\sqrt{5}}{3}) + h(1+\alpha^{m-n})&\leq& log3 + h(\alpha^{m-n}) + log2 \\&=& log6 + |m-n| h(\alpha)\\ &=& log6 + (n-m) \frac{\log \alpha}{2}.
\end{eqnarray*}
Next we prove that $|log\gamma3|<1$.
\begin{eqnarray*}
|\log\gamma3| = |\log\frac{\sqrt{5}}{3} (1+\alpha^{m-n})^{-1}|&=&|\log\frac{\sqrt{5}}{3}-log(1+\alpha^{m-n})|\\&<&|\log\frac{\sqrt{5}}{3}| + |\log(1+\alpha^{m-n})| 
\\&<&|\log\frac{\sqrt{5}}{3}| + |\log2| < 0.3 + 0.7=1 .
\end{eqnarray*} 
hence, we can take
\begin{equation*}
A_3=4 + (n-m) \log\alpha >max\{2h(\gamma_3) , |\log\gamma_3| , 0.16\}
\end{equation*}
Applying Matveev's theorem to (4.5), we get  
\begin{equation*}
\log |1-2^a \;\alpha^{-n}\; \frac{\sqrt{5}}{3} (1+\alpha^{m-n})^{-1}| > - C\: (1+\log n) \times 1.4 \times 0.5 \times( 4 +(n-m) \log\alpha) 
\end{equation*}
Where $C=1.4 \times 30^6 \times 3^{4.5} \times 2^2 \times (1+\log2) < 9.7\times 10^{11}$, using the fact that $(1+\log n) < 2\log n ,\text{for all}\; n \geq 3$ , we get 
\begin{equation*}
\log |1-2^a\; \alpha^{-n}\; \frac{\sqrt{5}}{3} (1+\alpha^{m-n})^{-1}| > -1.36 \times 10^{12} \times \log n \times (4 +(n-m) \log\alpha) 
\end{equation*}
By taking the logarithms in both sides in (4.5) and comparing the resulting inequality with the above relation, we get  
\begin{equation*}
log2 - n log\alpha > -1.36 \times 10^{12} \times log n \times (4 +(n-m) log\alpha)
\end{equation*}
The last inequality implies that 
\begin{equation}
n log \alpha < 1.37 \times 10^{12} \times log n \times (4 +(n-m) log\alpha)
\end{equation}
substituting by (4.4) in the above relation, we get   
\begin{equation*}
n log \alpha < 1.37 \times 10^{12} \times log n \times (4+3\times10^{12} \times log n)
\end{equation*}
solving the last inequality by using Mathematica, we get  
\begin{equation}
n<3.7 \times10^{28}
\end{equation}
\subsection{Reducing the bound on n} .
We need now to reduce the upper bound on n in inequality (4.7) to a size that can be easily handled, we use Lemma 2.1 several times to achieve the reduction.  
Consider 
\begin{equation}
 Z = a log2 - n log\alpha + log\frac{\sqrt{5}}{3}      
\end{equation}
Note that  we can write (4.3) as
\begin{equation*}
|1-e^z|<\frac{4}{\alpha^{n-m}}
\end{equation*}
Since we proved that the lift hand side of (4.3) is not equal zero,$|1-e^z|\neq0$, which implies that $Z \neq 0$. If $Z > 0$, we get 
\begin{equation}
0 < Z < e^z - 1 =|e^z - 1| < \frac{4}{\alpha^{n-m}}
\end{equation}
For $Z < 0$ , we suppose that\; $n-m \geq 20$, and since  $\frac{4}{\alpha^{n-m}} < \frac{1}{2}$ holds for  $(n-m) \geq20$, we get

\begin{equation*}
|1-e^z| =1-e^z < \frac{1}{2}
\end{equation*}
The last inequality implies that 
\begin{equation*}
e^{-z} < 2, \text{or}\; e^{|z|} < 2.
\end{equation*}
The last inequality implies that 
\begin{equation} 
0 < |Z| < e^{|z|} - 1 = e^{|z|} (1-e^z) = e^{|z|} |e^z -1| < \frac{8}{\alpha^{n-m}} 
\end{equation}

In both cases$(z>0, or \; z<0)$, (4.10) is true, replacing  Z in the last inequality by its formula (4.8) and dividing by $log\alpha$, we get 
\begin{equation}
0 < |a\frac{log 2}{log\alpha} - n + \frac{log\frac{\sqrt{5}}{3}}{log\alpha}| < \frac{17}{\alpha^{n-m}} 
\end{equation}

Taking
\begin{equation*}
\gamma = \frac{\log 2}{\log\alpha},\hspace{0.5cm} \mu = \frac{\log\frac{\sqrt{5}}{3}}{\log\alpha},\hspace{0.5cm} A = 17,\hspace{0.5cm} B = \alpha,\hspace{0.25cm}and\hspace{0.25cm} W = n-m .
\end{equation*}
It is clear that $\gamma$ is irrational number, from (4.7) we can take $M=3.7$×$10^{28}$ as an upper bound for n. For $q_{(67)} = 506642617666397667695263997821 > 6M = 2.22 \times 10^{29}$, where q is a denominator of a convergent of the continued fraction of $\gamma $. By the help of Mathematica, We find that $0.269087312907046<\epsilon=||\mu \:q|| -M\:||\gamma \: q|| < 0.269087312907048$. Applying Lemma 2.1 to inequality (4.11), we get 
\begin{equation}
 n-m<\frac{log(Aq/\epsilon)}{log B}<150.76
\end{equation}
inserting that upper bound for $n-m$ in (4.6), we get 
 \begin{equation}
 n < 8 \times 10^{15}
\end{equation}
Now we consider the other linear form (4.5)
\begin{equation*}
|1-2^a \alpha^{-n} \frac{\sqrt{5}}{3} (1+\alpha^{m-n})^{-1}| < \frac{2}{\alpha^n} 
\end{equation*}
Taking 
\begin{equation}
s = a log2 - n log\alpha + log \frac{\sqrt{5}}{3} (1+\alpha^{m-n})^{-1} 
\end{equation}
Equation (4.5) implies that 
\begin{equation} 
|1-e^s| < \frac{2}{\alpha^n} 
\end{equation}
Note that $s \neq 0$, since $|1-e^s|\neq0$. If $s > 0$, we get 
\begin{equation}
0 < s < e^s - 1 = |e^s - 1| < \frac{2}{\alpha^n}
\end{equation}
since $n > 200$,we get $\frac{2}{\alpha^n} < \frac{1}{2}$, which implies that For $s<0$  we have  
\begin{equation*}
|e^s -1| = 1- e^s < \frac{1}{2},\; \text{or} \;e^{|s|}=e^{-s} < 2.
\end{equation*}
The last inequality implies that 
\begin{equation}
0 < |s| < e^{|s|} - 1 = e^{|s|} (1-e^s) = e^{|s|} |e^s -1| < \frac{4}{\alpha^n}
\end{equation}
In both cases  ($s>0,s<0$), (2.15) is true. Replacing the value of s in the above inequality by its formula (4.14), and dividing by\; $\log\alpha$\;, we get 
\begin{equation}
0 < |a\frac{log 2}{log\alpha} - n + \frac{log \frac{\sqrt{5}}{3} (1+\alpha^{m-n})^{-1}}{log\alpha} |< \frac{9}{\alpha^{n}}
\end{equation}
Taking
\begin{equation*}
\gamma = \frac{log 2}{log\alpha},\hspace{0.5cm} \mu = \frac{log \frac{\sqrt{5}}{3} (1+\alpha^{m-n})^{-1}}{log\alpha},\hspace{0.5cm}w=n,\hspace{0.5cm} A = 9\hspace{0.25cm}\text{and}\hspace{0.25cm} B = \alpha.
\end{equation*}
and from (4.13) we can take $ M = 8 \times 10^{15} > n > a $ as an upper bound on a. Using the same $q_{(67)} = 506642617666397667695263997821 > 6M = 4.8 \times 10 ^{16} $, Mathematica finds that  in all cases for $n-m \in[1,150] $, The smallest epsilon is $0.0057323312747131<\epsilon= ||\mu \:q|| -M\:||\gamma \: q||<0.0057323312747133$, which belongs to $n-m=52$, applying Lemma 1.2 to (4.18), we get 
\begin{equation}
n<\frac{log(Aq/\epsilon)}{log B} <157.43
\end{equation}
In the last inequality we used the smallest epsilon because it produces the largest upper bound on n,(4.19) gives a contradiction to our assumption that $n \leq 200$, which shows that there are no solutions to (1.2) when $n > 200$.
\end{proof}
\makeatletter
\renewcommand{\@biblabel}[1]{[#1]\hfill}
\makeatother

\end{document}